\documentclass[11pt]{article}
\usepackage{amsmath,amssymb,graphicx}
\newcommand{\oC}{{\mathbb{C}}}
\newcommand{\oF}{{\mathbb{F}}}
\newcommand{\oN}{{\mathbb{N}}}
\newcommand{\oR}{{\mathbb{R}}}
\newcommand{\oZ}{{\mathbb{Z}}}
\newcommand{\AAA}{{\cal A}}
\newcommand{\GG}{{\cal G}}
\newcommand{\II}{{\cal I}}
\DeclareMathOperator{\rank}{rank}
\DeclareMathOperator{\height}{height}
\DeclareMathOperator{\sgn}{sgn}
\DeclareMathOperator{\GL}{GL}
\DeclareMathOperator{\tr}{tr}
\newcounter{bewering}
\newcommand{\prop}[2]{\refstepcounter{bewering}\vspace{4mm}\noindent{\bf Proposition \thebewering.}\label{#1}{\it #2}}
\newcommand{\kies}[2]{\mbox{${{#1}\choose{#2}}$}}
\newcommand{\pf}{\vspace{3mm}\noindent{\bf Proof.}\ }
\newcommand{\bx}{\hspace*{\fill} \hbox{\hskip 1pt \vrule width 4pt height 8pt depth 1.5pt \hskip 1pt}

\addvspace{4mm}}
\begin{document}
\begin{center}
{\large\bf CHARACTERIZING PARTITION FUNCTIONS OF THE EDGE-COLORING MODEL BY RANK GROWTH

}

\end{center}

\begin{center}
{\large
\hspace{10mm}
Alexander Schrijver\footnote{ University of Amsterdam and CWI, Amsterdam.
Mailing address: CWI, Science Park 123, 1098 XG Amsterdam,
The Netherlands.
Email: lex@cwi.nl.
The research leading to these results has received funding from the European Research Council
under the European Union's Seventh Framework Programme (FP7/2007-2013) / ERC grant agreement
n$\mbox{}^{\circ}$ 339109.}}

\end{center}

\noindent
{\small{\bf Abstract.}
We characterize which graph invariants are partition functions
of an edge-coloring model over $\oC$, in terms of the rank growth of
associated `connection matrices'.

}

\section{Introduction}

Let $\GG$ denote the collection of all undirected graphs, two of them
being the same if they are isomorphic.
In this paper, all graphs are finite and may have loops and multiple edges.
Let $k\in\oN$ and let $\oF$ be a commutative ring.
Call any function $y:\oN^k\to\oF$ a
({\em $k$-color}) {\em edge-coloring model} ({\em over $\oF$}).
In the case where $y$ is symmetric under the action of $S_k$, an edge-coloring
model is called a `vertex model' by de la Harpe and Jones [5], where
colors are called `states'.
The more general model was considered in the context of Holant functions by
L.G. Valiant (cf.\ [1]).

The {\em partition function} of an edge-coloring model $y$ is the function $p_y:\GG\to\oF$
defined for any graph $G=(V,E)$ by
\begin{align*}
p_y(G):=\sum_{\kappa:E\to[k]}\prod_{v\in V}y_{\kappa(\delta(v))}.
\end{align*}
Here $\delta(v)$ is the set of edges incident with $v$.
Then $\kappa(\delta(v))$ is a multisubset of $[k]$,
which we identify with its incidence vector in $\oN^k$.
Moreover, we use $\oN=\{0,1,2,\ldots\}$ and for $n\in\oN$,
\begin{align*}
[n]:=\{1,\ldots,n\}.
\end{align*}

We can visualize $\kappa$ as a coloring of the edges of $G$ and
$\kappa(\delta(v))$ as the multiset of colors `seen' from $v$.
The edge-coloring model was considered by
de la Harpe and Jones [5] as a physical model, where vertices serve
as particles, edges as interactions between particles,
and colors as states or energy levels.
It extends the Ising-Potts model.
Several graph parameters are partition functions of some edge-coloring model,
like the number of matchings.
There exist real-valued graph parameters that are partition functions
of an edge-coloring model over $\oC$, but not over $\oR$.
(A simple example is $(-1)^{|E(G)|}$.)

In this paper, we characterize which functions $f:\GG\to\oC$
are the partition function of an edge-coloring model over $\oC$.
The characterization differs from an earlier characterization
given in [3]
(which our present characterization uses) in that it is
based on the rank growth of associated `connection matrices'.

To describe it, we need the notion of a $k$-fragment.
For $k\in\oN$, a {\em $k$-fragment} is an undirected graph $G=(V,E)$
together with an injective `label' function $\lambda:[k]\to V$,
where $\lambda(i)$ is a vertex of degree 1, for each $i\in[k]$.
(You may alternatively view these degree-1 vertices
as ends of `half-edges', or rather of `edge pieces', as both ends of
an edge might be labeled.)

If $G$ and $H$ are $k$-fragments, the graph $G\cdot H$
is obtained from the disjoint union of $G$ and $H$ by identifying equally
labeled vertices and by ignoring each of the $k$ identified points
as vertex, joining its two incident edges into one edge.
\begin{align*}
\raisebox{-.45\height}{\scalebox{0.2}{\includegraphics{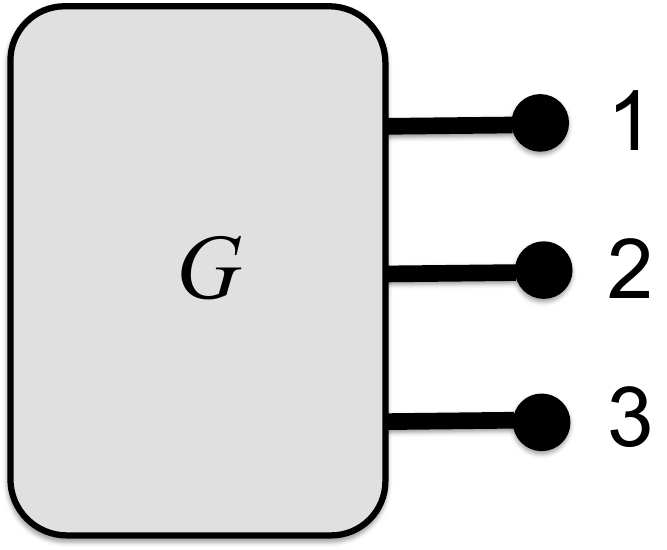}}}
~~~~
\cdot
~~~~
\raisebox{-.45\height}{\scalebox{0.2}{\includegraphics{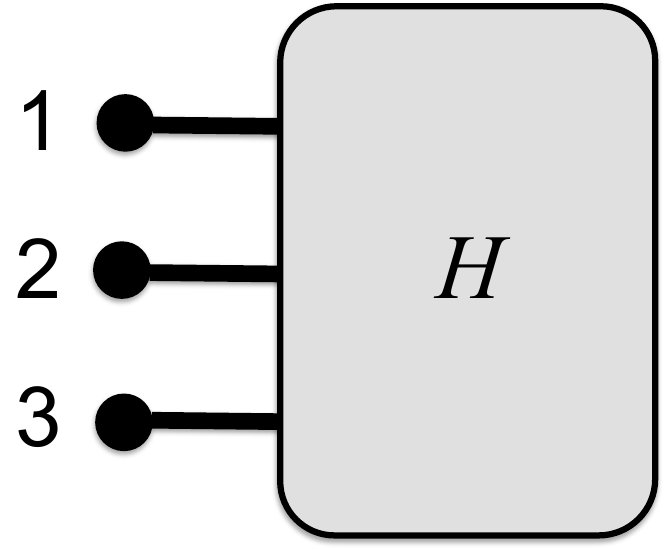}}}
~~~~
=
~~~~
\raisebox{-.45\height}{\scalebox{0.2}{\includegraphics{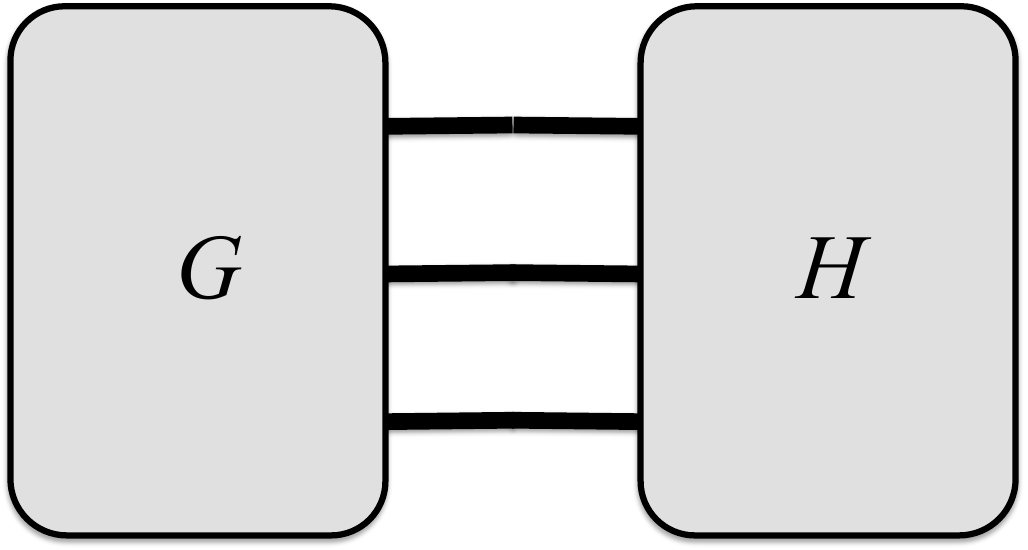}}}
\end{align*}
\begin{center}
The multiplication $G\cdot H$.
\end{center}
(A good way to imagine this is to see a graph as a topological $1$-complex.)
Note that it requires (as in [8])
that we also should consider the `vertexless
loop' as possible edge of a graph, as we may create it in $G\cdot H$.
We denote this vertexless loop by $\bigcirc$.
Observe that if $y$ is an edge-coloring model over $\oC$ with $n$ colors, then $p_y(\bigcirc)=n$.

Let $\GG_k$ denote the collections of $k$-fragments.
For any $f:\GG\to\oC$ and $k\in\oN$, the {\em $k$th connection matrix}
is the $\GG_k\times\GG_k$ matrix $C_{f,k}$ defined by
\begin{align*}
(C_{f,k})_{G,H}:=f(G\cdot H)
\end{align*}
for $G,H\in\GG_k$.

Now we can formulate our characterization:

\medskip
\noindent
{\bf Theorem.}{\it
A function $f:\GG\to\oC$
is the partition function of an edge-coloring model over $\oC$
if and only if $f(\emptyset)=1$, $f(\bigcirc)\in\oR$, and
\begin{align*}
\rank(C_{f,k})\leq f(\bigcirc)^k
\end{align*}
for each $k\in\oN$.
}

\medskip
Let us relate this to Szegedy's theorem [8], which
characterizes the partition functions of edge-coloring models over $\oR$.
Call a function $f:\GG\to\oC$ {\em multiplicative} if $f(\emptyset)=1$
and $f(G\dot{\cup} H)=f(G)f(H)$ for all graphs $G$ and $H$, where $G\dot{\cup} H$ denotes
the disjoint union of $G$ and $H$.
Then Szegedy's theorem reads:
\begin{align*}
\parbox{4in}{
A function $f:\GG\to\oR$ is the partition function of
an edge-coloring model over $\oR$
if and only if $f$ is multiplicative and
$C_{f,k}$ is positive semidefinite for each $k$.
}
\end{align*}
For related results for the `spin model' see [4] and [7].

The proof of our theorem is based on some elementary
results from the representation theory
of the symmetric group, and on the following alternative
characterization of partition
functions of edge-coloring models given in
[3], which uses the Nullstellensatz
and the First and Second Fundamental Theorems of Invariant theory for $O(n)$.

For any graph $G=(V,E)$, any $U\subseteq V$, and any $s:U\to V$, define
\begin{align*}
E_s:=\{us(u)\mid u\in U\}\text{ and }G_s:=(V,E\cup E_s),
\end{align*}
where $us(u)$ denotes an edge connecting $u$ and $s(u)$
(adding multiple edges if $E_s$ intersects $E$).
Let $S_U$ be the group of permutations of $U$.
Then ([3]):
\begin{align}\label{9no10c}
\parbox[t]{4in}{A function $f:\GG\to\oC$ is the partition function of some
$k$-color edge-coloring model over $\oC$
if and only if $f$ is multiplicative
and for each graph $G=(V,E)$, each $U\subseteq V$ with $|U|=k+1$, and
each $s:U\to V$:
$\displaystyle
\sum_{\pi\in S_U}\sgn(\pi)f(G_{s\circ\pi})=0$.
}
\end{align}

\section{Some results on the symmetric group}

In the proof of our theorem we will need Proposition \ref{26jl11c} below, which we prove
in a number of steps.
(The result might be known, and must not be a difficult exercise for those familiar
with the representation theory of the symmetric group, but we did
not find an explicit reference.)

We recall a few standard results from the representation theory of
the symmetric group $S_n$ (cf.\ James and Kerber [6]).
Basis is the one-to-one relation between the partitions $\lambda$ of $n$ and the
irreducible representations $r_{\lambda}$ of $S_n$.
Here a {\em partition} $\lambda$ of $n$ is a finite nonincreasing sequence
$(\lambda_1,\ldots,\lambda_t)$ of positive integers with sum $n$.
One writes $\lambda\vdash n$ if $\lambda$ is a partition of $n$.
The number $t$ of terms of $\lambda$ is called the {\em height} of
$\lambda$, denoted by $\height(\lambda)$.
Denote by $f^{\lambda}$ the degree of representation $r_{\lambda}$
(that is, the dimension of the representation space of $f^{\lambda}$), and by
$\chi_{\lambda}$ the character of $r_{\lambda}$.

For any $\lambda\vdash n$, the {\em Young shape} $Y_{\lambda}$ of $\lambda=(\lambda_1,\ldots,\lambda_t)$ is
the following subset of $\oN^2$:
\begin{align*}
Y_{\lambda}:=\{(i,j)\mid i\in[t], j\in[\lambda_i]\}.
\end{align*}
For any $\pi\in S_n$, let $o(\pi)$ denote the number of orbits of $\pi$.

\prop{26jl11e}{
For any $n\in\oN$, $\lambda\vdash n$, and $d\in\oC$:
\begin{align}\label{26jl11f}
\sum_{\pi\in S_n}\chi_{\lambda}(\pi)d^{o(\pi)}
=
f^{\lambda}\prod_{(i,j)\in Y_{\lambda}}(d+j-i).
\end{align}
}

\pf
As both sides of \eqref{26jl11f} are polynomials in $d$, we can assume
that $d\in\oN$.
Consider the representation $r$ of $S_n$ on $(\oC^d)^{\otimes n}$ induced
by
\begin{align*}
\pi\cdot(x_1\otimes\cdots\otimes x_n)=x_{\pi(1)}\otimes\cdots\otimes x_{\pi(n)}
\end{align*}
for $x_1,\ldots,x_n\in \oC^d$.
Note that the character $\chi$ of $r$ satisfies $\chi(\pi)=d^{o(\pi)}$ for
each $\pi\in S_n$, since, fixing a basis $e_1,\ldots,e_d$ of $\oC^d$, it
is equal to the number of $z=e_{i_1}\otimes\cdots\otimes e_{i_n}$ with
$\pi\cdot z=z$; that is,
it is equal to the number of $(i_1,\ldots,i_n)\in[d]^n$
with $(\pi(i_1),\ldots,\pi(i_n))=(i_1,\ldots,i_n)$.
So the requirement is that $i_j=i_k$ whenever
$j$ and $k$ belong to the same orbit of $\pi$.
Therefore, this number is equal to $d^{o(\pi)}$.

For any $\alpha\vdash n$, let
$\mu_{\alpha}$ be the multiplicity of $r_{\alpha}$ in $r$.
Then, using Schur-Weyl duality,
\begin{align*}
&\sum_{\pi\in S_n}
\chi_{\lambda}(\pi)
d^{o(\pi)}
=
\sum_{\pi\in S_n}
\chi_{\lambda}(\pi)
\overline{\chi(\pi)}
=
\sum_{\pi\in S_n}
\chi_{\lambda}(\pi)
\sum_{\alpha\vdash n}
\mu_{\alpha}\overline{\chi_{\alpha}(\pi)}\\
&=
\sum_{\alpha\vdash n}
\mu_{\alpha}
\sum_{\pi\in S_n}
\chi_{\lambda}(\pi)
\overline{\chi_{\alpha}(\pi)}
=
\sum_{\alpha\vdash n}\mu_{\alpha}n!\delta_{\lambda,\alpha}
=
n!\mu_{\lambda}\\
&=
f^{\lambda}\prod_{(i,j)\in Y_{\lambda}}(d+j-i).
\end{align*}
The last equality follows from the fact that $\mu_{\lambda}$ is
equal to the degree of the irreducible
representation of $\GL(d,\oC)$ corresponding to $\lambda$
(by Schur-Weyl duality)
and that this degree is equal to the last expression divided by $n!$
(cf.\ [2] eq.\ 9.28).
This shows \eqref{26jl11f}.
\bx

For any $n\in\oN$ and $d\in\oC$, let $M_n(d)$ be the $S_n\times S_n$
matrix with
\begin{align*}
(M_n(d))_{\rho,\sigma}:=d^{o(\rho\sigma^{-1})}
\end{align*}
for $\rho,\sigma\in S_n$.

\prop{26jl11d}{
For any $n\in\oN$ and $d\in\oC$:
\begin{align}\label{26jl11g}
\rank(M_n(d))=
\begin{cases}
n!&\text{ if $d\not\in\oZ$,}\\
\sum((f^{\lambda})^2\mid\lambda\vdash n, \height(\lambda)\leq |d|)&\text{ if $d\in\oZ$.}
\end{cases}
\end{align}
}

\pf
First, we have, for any $d\in\oC$,
\begin{align}\label{27jl11a}
\rank(M_n(-d))=\rank(M_n(d)).
\end{align}
Indeed, note that $M_n(-d)=(-1)^n\Delta_{\sgn}M_n(d)\Delta_{\sgn}$,
where $\Delta_{\sgn}$ is the $S_n\times S_n$ diagonal matrix
with $(\Delta_{\sgn})_{\pi,\pi}=\sgn(\pi)$ for $\pi\in S_n$.
(This because $\sgn(\pi)=(-1)^{n-o(\pi)}$ for all $\pi$,
hence $(-1)^{o(\rho\sigma^{-1})}
=
(-1)^n\sgn(\rho\sigma^{-1})
=
(-1)^n\sgn(\rho)\sgn(\sigma)$.)
This gives \eqref{27jl11a}.

Let $R$ be the regular representation of $S_n$.
So, for any $\pi\in S_n$, $R(\pi)$ is the $S_n\times S_n$ matrix
with
\begin{align*}
R(\pi)_{\rho,\sigma}=
\begin{cases}
1&\text{ if $\rho=\pi\sigma$,}\\
0&\text{ otherwise}
\end{cases}
\end{align*}
for $\rho,\sigma\in S_n$.
Then
\begin{align*}
M_n(d)=\sum_{\pi\in S_n}d^{o(\pi)}R(\pi).
\end{align*}
Now $M_n(d)$ commutes with each $R(\pi)$, and belongs
to the group algebra of $S_n$.
So $\rank(M_n(d))$ is equal to the sum of $(f^{\lambda})^2$
taken over those $\lambda$ such that $M_n(d)$ has nonzero
trace in representation $r_{\lambda}$.
That is, such that
$\sum_{\pi\in S_n}\chi_{\lambda}(\pi)d^{o(\pi)}\neq 0$.
So
\begin{align*}
&\rank(M_n(d))=
\sum((f^{\lambda})^2
\mid
\lambda\vdash n,
\sum_{\pi\in S_n}\chi_{\lambda}(\pi)d^{o(\pi)}\neq 0)\\
&=
\sum((f^{\lambda})^2
\mid
\lambda\vdash n,
d\not\in\{i-j\mid (i,j)\in Y_{\lambda}\}).
\end{align*}
The last equality follows from \eqref{26jl11f}.

Now if $d\not\in\oZ$, then for all $\lambda\vdash n$:
$d\neq i-j$ for all $(i,j)\in Y_{\lambda}$.
So $\rank(M_{n,d})=n!$.
If $d\in\oZ$, then by \eqref{27jl11a} we can assume that $d$ is nonnegative,
that is, $d\in\oN$.
Then for each $\lambda\vdash n$:
$d\not\in \{i-j\mid (i,j)\in Y_{\lambda}\}$
if and only if $\height(\lambda)\leq d$.
This proves \eqref{26jl11g}.
\bx

\prop{26jl11c}{
For any $d\in\oC$:
\begin{align}\label{26jl11h}
\sup_{n\in\oN}
\left(\rank(M_n(d))\right)^{1/n}=
\begin{cases}
\infty&\text{ if $d\not\in\oZ$,}\\
d^2&\text{ if $d\in\oZ$.}
\end{cases}
\end{align}
}

\pf
If $d\not\in\oZ$, the result follows directly from \eqref{26jl11g},
as $\sup_nn!^{1/n}=\infty$.

If $d\in\oZ$, then again by \eqref{27jl11a} we can assume that $d$ is nonnegative,
that is, $d\in\oN$.
Then $\rank(M_n(d))\leq (d^2)^n$.
Indeed, let $\chi$ be the character of the natural representation $r$
of $S_n$ on $(\oC^d)^{\otimes n}$.
Then $d^{o(\pi)}=\chi(\pi)$ for all $\pi\in S_n$.
Hence $d^{o(\rho\sigma^{-1})}=\chi(\rho\sigma^{-1})$.
So $d^{o(\rho\sigma^{-1})}$ is the trace of the product of the
$d^n\times d^n$ matrices $r(\rho)$ and $r(\sigma^{-1})$.
Hence $\rank(M_n(d))\leq (d^n)^2=(d^2)^n$.
This proves $\leq$ in \eqref{26jl11h}.

To prove the reverse inequality, consider for any $m\in\oN$,
the partition $\lambda_m=(m,\ldots,m)$ of $n:=dm$, with $\height(\lambda_m)=d$.
By the hook formula,
\begin{align*}
f^{\lambda_m}=n!/\prod_{i=1}^d\prod_{j=1}^m(i+j-1)
=
\frac{(dm)!0!1!\cdots(d-1)!}{m!(m+1)!\cdots(m+d-1)!}
=
\frac{(dm)!}{m!^dp(m)},
\end{align*}
where (fixing $d$) $p(m)$ is a polynomial in $m$
(namely $p(m)=\prod_{i=0}^{d-1}\kies{m+i}{i}$).
So, by Stirling's formula, $\lim_{m\to\infty}(f^{\lambda_m})^{1/dm}=d$.
By \eqref{26jl11g}, we have for each $m$, since $\lambda_m\vdash dm$ and
$\height(\lambda_m)=d$,
\begin{align*}
\rank(M_{dm}(d))\geq (f^{\lambda_m})^2.
\end{align*}
This gives the required inequality.
\bx

\section{Proof of the theorem}

Necessity being easy, we show sufficiency.
As $f(\emptyset)=1$ and $\rank(C_{f,0})\leq f(\bigcirc)^0=1$, we know that
$f$ is multiplicative.
Moreover, as $\rank(C_{f,1})\leq f(\bigcirc)$, we know $f(\bigcirc)\geq 0$.

We develop some straightforward algebra.
Let $k\in\oN$.
For $G,H\in\GG_{2k}$, define the product $GH$ as the $2k$-fragment obtained
from the disjoint union of $G$ and $H$ by identifying vertex labeled $k+i$
in $G$ with vertex labeled $i$ in $H$, and ignoring this vertex as vertex
(for $i=1,\ldots,k$);
the vertices of $G$ labeled $1,\ldots,k$ and those of $H$ labeled $k+1,\ldots,2k$
make $GH$ to a $2k$-fragment again.
\begin{align*}
\raisebox{-.45\height}{\scalebox{0.19}{\includegraphics{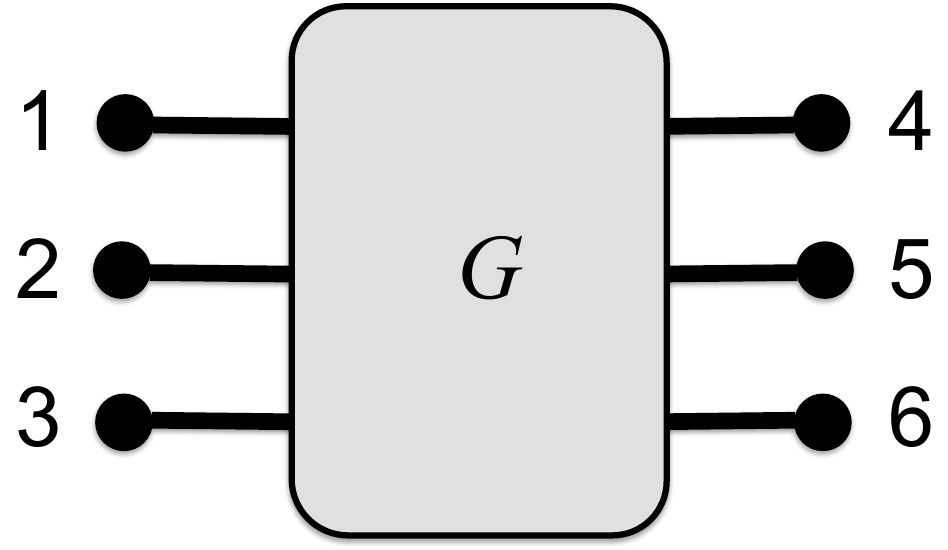}}}
~~
\raisebox{-.45\height}{\scalebox{0.19}{\includegraphics{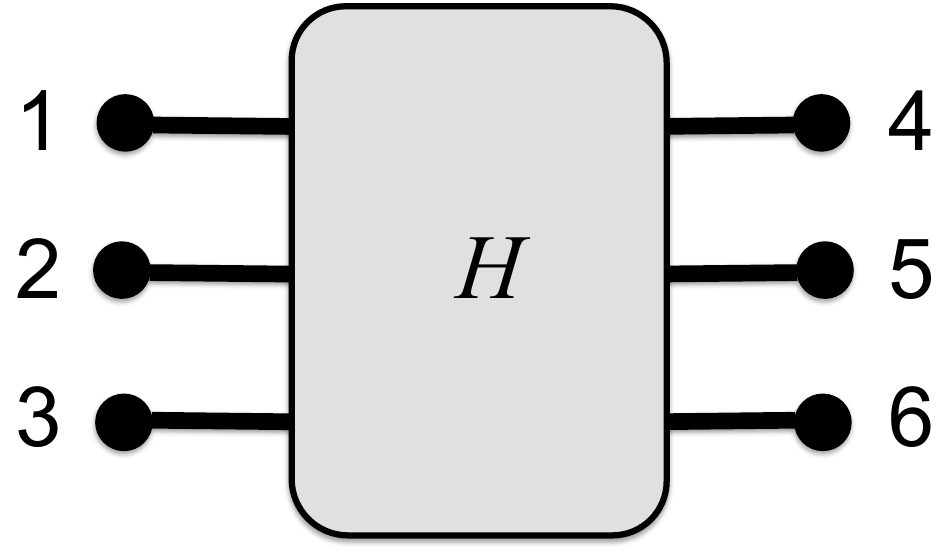}}}
~~
=
~~
\raisebox{-.45\height}{\scalebox{0.19}{\includegraphics{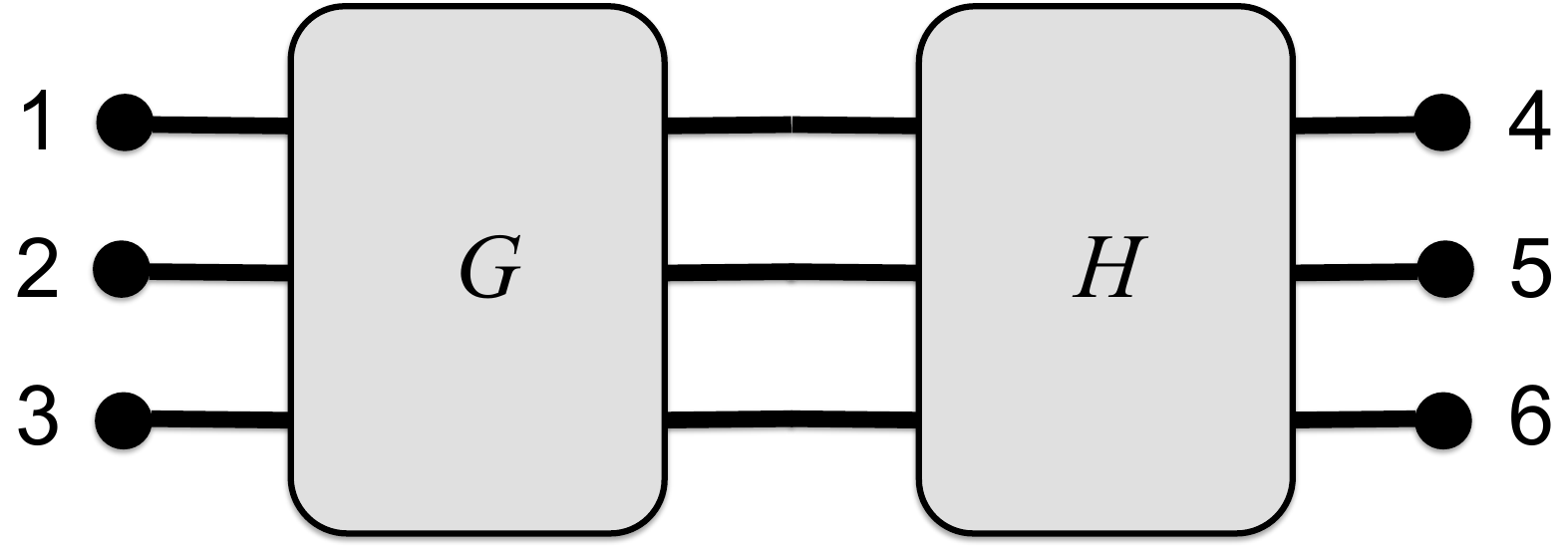}}}
\end{align*}
\begin{center}
The multiplication $GH$.
\end{center}

Geometrically, one may imagine that the $2k$-fragments have the labels $1,\ldots,k$
vertically at the left and the labels $k+1,\ldots,2k$ vertically at the right.
Then $GH$ arises by drawing $G$ at the left from $H$ and connecting the right-side
labels of $G$ with the left-side labels of $H$, in order.

Clearly, this product is associative.
Moreover, there is a unit, denoted by ${\mathbf 1}_k$, consisting of $k$
disjoint edges $e_1,\ldots,e_k$, where the ends of $e_i$ are labeled $i$
and $k+i$ ($i=1,\ldots,k$).

Let $\oC\GG_{2k}$ be the collection of formal $\oC$-linear combinations
of elements of $G_{2k}$.
Extend the products $G\cdot H$ and $GH$ bilinearly to $\oC\GG_{2k}$.
The latter product makes $\oC\GG_{2k}$ to a $\oC$-algebra.

Let $\II_{2k}$ be the null space of the matrix $C_{f,2k}$; that is, it consists of all
$\gamma\in\oC\GG_{2k}$ with $f(\gamma\cdot H)=0$ for all $H\in\GG_{2k}$.
Then $\II_{2k}$ is an ideal in the algebra $\oC\GG_{2k}$,
and the quotient
\begin{align*}
\AAA_k:=\oC\GG_{2k}/\II_{2k}
\end{align*}
is an algebra of dimension $\rank(C_{f,2k})$.
We will indicate elements of $\AAA_k$ by representatives in $\oC\GG_{2k}$.

Define the `trace-like' function $\tau:\AAA_k\to\oC$ by
\begin{align*}
\tau(x):=f(x\cdot {\mathbf 1}_k).
\end{align*}
Then $\tau(xy)=\tau(yx)$ for all $x,y\in\AAA_k$ and $\tau({\mathbf 1}_k)=f(\bigcirc)^k$.
Note also that if $G,H\in\GG_{2k}$, then $f(G\cdot H)=\tau(G\widetilde H)$,
where $\widetilde H$ arises from $H$ by exchanging labels
$i$ and $k+i$ for each $i=1,\ldots,k$.
Extending this linearly to $\oC\GG_{2k}$, we know that for each
$x\in\AAA_k$:
\begin{align}\label{18me15c}
\parbox[t]{3.5in}{
if $f(x\cdot z)\neq 0$ for some $z\in\AAA_k$, then
$\tau(xy)\neq 0$ for some $y\in\AAA_k$.
}
\end{align}

We will first show that $\AAA_k$ is semisimple.
To this end, for $G\in\GG_{2k}$ and $H\in\GG_{2l}$, we need a product $G\sqcup H$ that can
be intuitively described as: put $G$ above $H$, and renumber the labels at the left hand side
to $1,\ldots,k+l$ (in order) and renumber the labels at the right hand side to
$k+l+1,\ldots,2(k+l)$ (in order).
\begin{align*}
\raisebox{-.45\height}{\scalebox{0.19}{\includegraphics{gremrk_G2.pdf}}}
~~
\sqcup
~~
\raisebox{-.45\height}{\scalebox{0.19}{\includegraphics{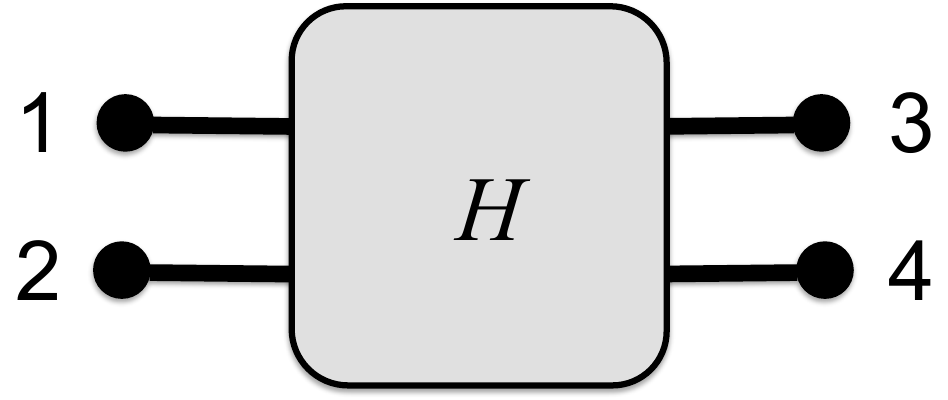}}}
~~
=
~~
\raisebox{-.45\height}{\scalebox{0.19}{\includegraphics{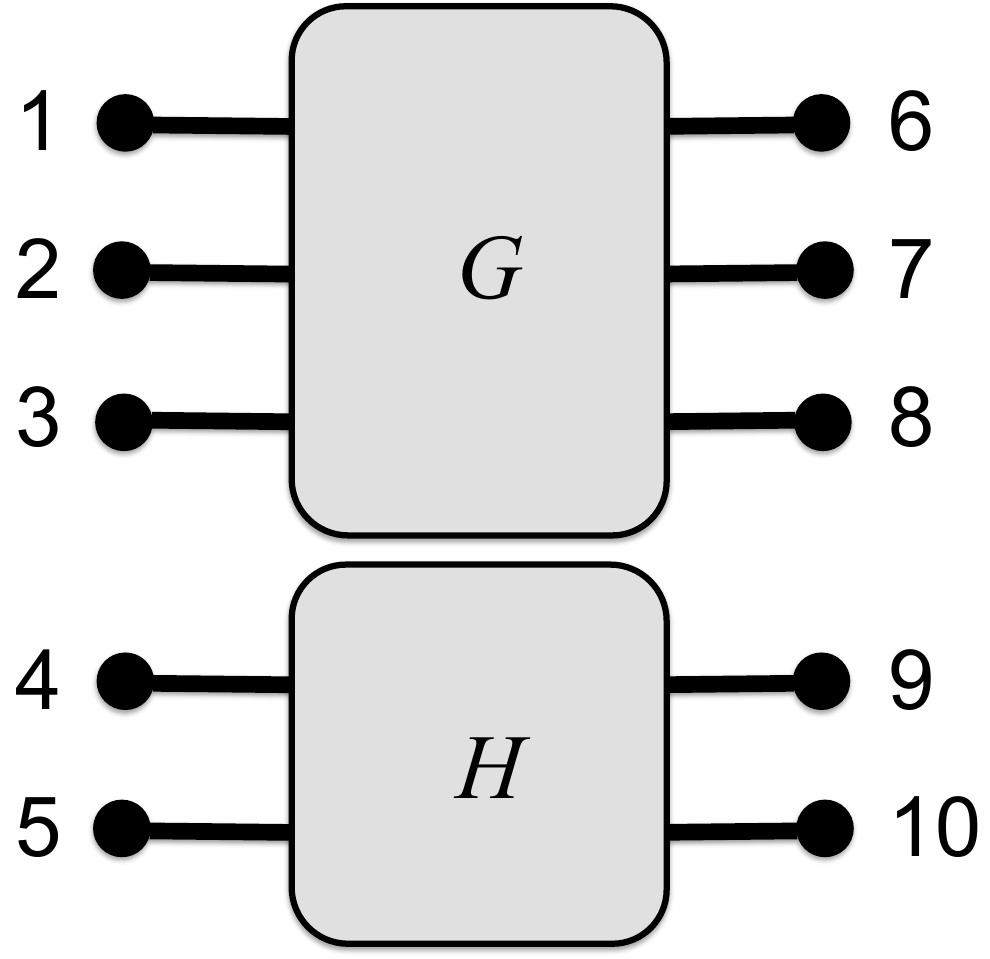}}}
\end{align*}
\begin{center}
The multiplication $G\sqcup H$.
\end{center}

More precisely stated, $G\sqcup H$ is the $2(k+l)$-fragment
obtained from the disjoint union of $G$ and $H$ by adding $l$ to the labels $k+1,\ldots,2k$
of $G$, adding $k$ to the labels $1,\ldots,l$ in $H$, and by adding $2k+l$ to the labels
$l+1,\ldots,2l$ in $H$.
This product $\sqcup$ is associative and extends bilinearly to
$\oC\GG_{2k}\times\oC\GG_{2l}\to\oC\GG_{2(k+l)}$.
Thus for $x\in\oC\GG_{2k}$, the $m$-th power $x^{\sqcup m}$ is well-defined.

Consider any $k,m\in\oN$.
For $\pi\in S_m$, let $P_{k,\pi}$ be the $2km$-fragment consisting of $km$ disjoint
edges $e_{i,j}$ for $i=1,\ldots,m$ and $j=1,\ldots,k$, where $e_{i,j}$ connects
the vertices labeled $(i-1)k+j$ and $km+(\pi(i)-1)k+j$.
\begin{align*}
\raisebox{-.45\height}{\scalebox{0.19}{\includegraphics{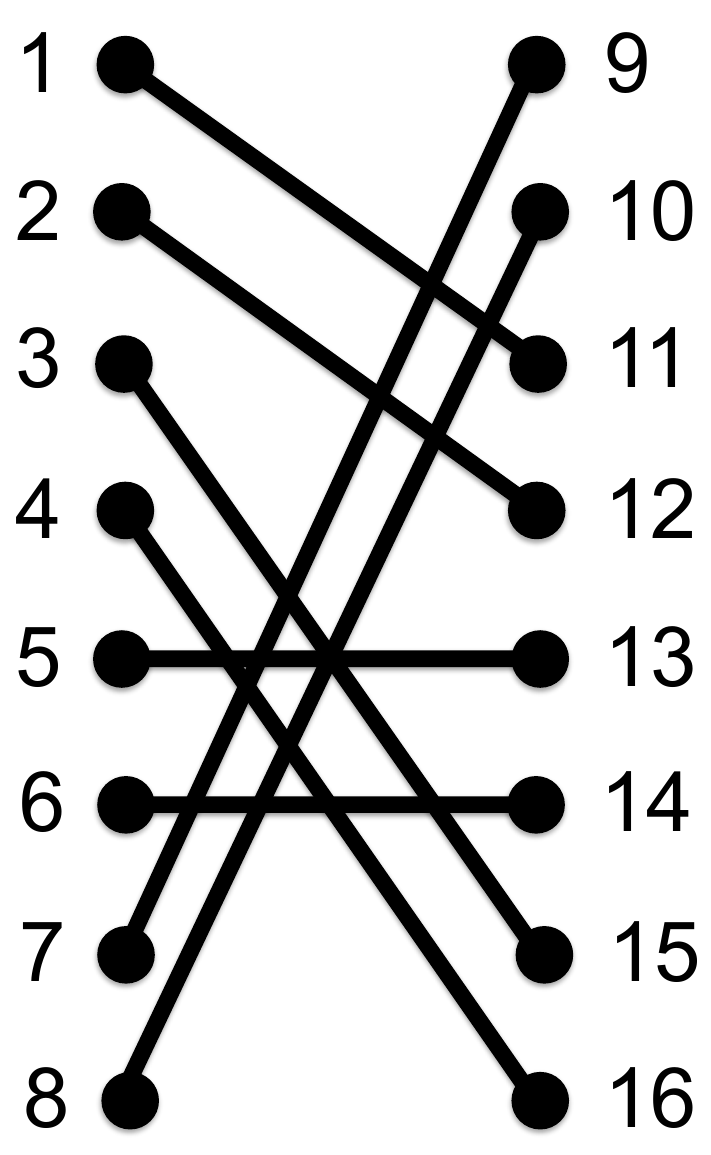}}}
\end{align*}
\begin{center}
The $16$-fragment $P_{2,\pi}$ with $\pi=(124)\in S_4$.
\end{center}

Then for any $\rho,\sigma\in S_m$ one has
\begin{align*}
f(x^{\sqcup m}\cdot P_{k,\pi})
=
\prod_c\tau(x^{|c|}),
\end{align*}
where $c$ ranges over the orbits of $\pi$.
Hence for any $\rho,\sigma\in S_m$ one has
\begin{align}\label{5se11a}
f(x^{\sqcup m}P_{k,\rho}\cdot P_{k,\sigma})
=
f(x^{\sqcup m}\cdot P_{k,\rho^{-1}}P_{k,\sigma})
=
f(x^{\sqcup m}\cdot P_{k,\rho^{-1}\sigma})
=
\prod_c\tau(x^{|c|}),
\end{align}
where $c$ ranges over the orbits of $\rho^{-1}\sigma$.
We will use that, for each $x\in\oC\GG_{2k}$, the $S_m\times S_m$ matrix
$(f(x^{\sqcup m}P_{k,\rho}\cdot P_{k,\sigma})_{\rho,\sigma\in S_m}$ has rank at most
$\rank(C_{f,2k})$ (since $x^{\sqcup m}P_{k,\rho}$ belongs to $\oC\GG_{2k}$ for each $\rho$).

\prop{5se11b}{
If $x$ is a nilpotent element of $\AAA_k$, then $\tau(x)=0$.
}

\pf
Suppose $\tau(x)\neq 0$ and $x$ is nilpotent.
Then there is a largest $t$ with $\tau(x^t)\neq 0$.
Let $y:=x^t$.
So $\tau(y)\neq 0$ and $\tau(y^s)=0$ for each $s\geq 2$.
By scaling we can assume that $\tau(y)=1$.

Choose $m$ such that $m!>f(\bigcirc)^{2km}$.
By \eqref{5se11a} we have, for any $\rho,\sigma\in S_m$,
\begin{align*}
f(y^{\sqcup m}P_{k,\rho}\cdot P_{k,\sigma})=\delta_{\rho,\sigma},
\end{align*}
since \eqref{5se11a} is 0 if $\rho^{-1}\sigma$ has an orbit $c$ with $|c|>1$, i.e., $\rho^{-1}\sigma$ is not the identity.

So $\rank(C_{f,2km})\geq m!$, contradicting the fact that
$\rank(C_{f,2km})\linebreak[4]\leq f(\bigcirc)^{2km}<m!$.
\bx

The following is a direct consequence of Proposition \ref{5se11b}:

\prop{18me15a}{
$\AAA_k$ is semisimple.
}

\pf
As $\AAA_k$ is finite-dimensional, it suffices to
show that for each nonzero element $x$ of $\AAA_k$ there is a
$y$ with $xy$ not nilpotent.
As $x\not\in\II_{2k}$, we know that $f(x\cdot z)\neq 0$ for some
$z\in \AAA_k$.
So by \eqref{18me15c}, $\tau(xy)\neq 0$ for some $y\in\AAA_k$, and hence,
by Proposition \ref{5se11b}, $xy$ is not nilpotent.
\bx

\prop{26jl11j}{
If $x$ is a nonzero idempotent in $\AAA_k$, then $\tau(x)$ is
a positive integer.
}

\pf
Let $x$ be any idempotent.
Then for each $m\in\oN$ and $\rho,\sigma\in S_m$, by \eqref{5se11a}:
\begin{align*}
f(x^{\sqcup m}P_{k,\rho}\cdot P_{k,\sigma})=\tau(x)^{o(\rho\sigma^{-1})}.
\end{align*}
So for each $m$:
\begin{align*}
\rank(M_{m}(\tau(x)))\leq\rank(C_{f,2km})\leq f(\bigcirc)^{2km}.
\end{align*}
Hence
\begin{align*}
\sup_m(\rank(M_m(\tau(x))))^{1/m}\leq f(\bigcirc)^{2k}.
\end{align*}
By Proposition \ref{26jl11c} this implies $\tau(x)\in\oZ$ and $\tau(x)\leq f(\bigcirc)^k$.
As ${\mathbf 1}_k-x$ also is an idempotent in $\oC\GG_{2k}$
and as $\tau({\mathbf 1}_k)=f(\bigcirc)^k$,
we have
\begin{align*}
f(\bigcirc)^k\geq\tau({\mathbf 1}_k-x)=f(\bigcirc)^k-\tau(x).
\end{align*}
So $\tau(x)\geq 0$.

Suppose finally that $x$ is nonzero while $\tau(x)=0$.
As $\tau(y)\geq 0$ for each idempotent $y$,
we may assume that $x$ is a minimal nonzero idempotent.
Let $J$ be the two-sided ideal generated by $x$.
As $\AAA_k$ is semisimple, $J\cong\oC^{m\times m}$ for some $m$.
As $\tau$ is linear, there exists an $a\in J$ such that
$\tau(z)=\tr(za)$ for each $z\in J$.
As $\tau(z)=0$ for each nilpotent $z$, we know that $a$ is a diagonal matrix.
As $\tau(yz)=\tau(zy)$ for all $y,z\in J$, $a$ is in fact equal to
a scalar multiple of the identity matrix.

As $x\neq 0$, $f(x\cdot z)\neq 0$ for some $z\in\AAA_k$.
So by \eqref{18me15c}, $\tau(xy)\neq 0$ for some $y$.
Hence $a\neq 0$, and so $\tau(x)\neq 0$, contradicting our assumption.
\bx

As ${\mathbf 1}_1$ is an idempotent, we know that
$\tau({\mathbf 1}_1)$ is a nonnegative integer, say $n$.
So $f(\bigcirc)=n$.
Let $k:=n+1$.
For $\pi\in S_k$ let $r_{\pi}$ be the $2k$-fragment
consisting of $k$ disjoint edges $e_1,\ldots,e_k$, where
the ends of $e_i$ are labeled $i$ and $k+\pi(i)$,
for $i=1,\ldots,k$.
(In fact, $r_{\pi}=P_{1,\pi}$ as defined above.)
We define the following element $q$ of $\oC\GG_{2k}$:
\begin{align*}
q:=\sum_{\pi\in S_k}\sgn(\pi)r_{\pi}.
\end{align*}
By \eqref{9no10c} it suffices to show that $q\in\II_{2k}$,
that is, $q=0$ in $\AAA_k$.

Now $k!^{-1}q$ is an idempotent in $\oC\GG_{2k}$.
Moreover,
\begin{align*}
\tau(q)
=
\sum_{\pi\in S_k}\sgn(\pi)n^{o(\pi)}
=
\sum_{\pi\in S_k}\sgn(\pi)\sum_{\varphi:[k]\to[n]\atop\varphi\circ\pi=\varphi}1
=
\sum_{\varphi:[k]\to[n]}
\sum_{\pi\in S_k\atop\varphi\circ\pi=\varphi}\sgn(\pi)
=
0,
\end{align*}
since no $\varphi:[k]\to[n]$ is injective.
So by Proposition \ref{26jl11j}, $q=0$ in $\AAA_k$, as required.
This finishes the proof of the theorem.

\medskip
\noindent
{\bf Acknowledgements.}
I am grateful to the referees for very helpful suggestions as to the
presentation of this paper.

\section*{References}
{\small
\begin{itemize}{}{
\setlength{\labelwidth}{4mm}
\setlength{\parsep}{0mm}
\setlength{\itemsep}{1mm}
\setlength{\leftmargin}{5mm}
\setlength{\labelsep}{1mm}
}
\item[\mbox{\rm[1]}] J.-Y. Cai, P. Lu, M. Xia,
Holographic algorithms by Fibonacci gates,
{\em Linear Algebra and its Applications} 438 (2013) 690--707.

\item[\mbox{\rm[2]}] P. Cvitanovi\'c,
{\em Group Theory: Birdtracks, Lie's, and Exceptional Groups},
Princeton University Press, Princeton, 2008.

\item[\mbox{\rm[3]}] J. Draisma, D. Gijswijt, L. Lov\'asz, G. Regts, A. Schrijver,
Characterizing partition functions of the vertex model,
{\em Journal of Algebra} 350 (2012) 197--206.

\item[\mbox{\rm[4]}] M.H. Freedman, L. Lov\'asz, A. Schrijver,
Reflection positivity, rank connectivity, and homomorphisms of graphs,
{\em Journal of the American Mathematical Society} 20 (2007) 37--51.

\item[\mbox{\rm[5]}] P. de la Harpe, V.F.R. Jones,
Graph invariants related to statistical mechanical models
examples and problems,
{\em Journal of Combinatorial Theory, Series B} 57 (1993) 207--227.

\item[\mbox{\rm[6]}] G. James, A. Kerber,
{\em The Representation Theory of the Symmetric Group},
Addison-Wesley, Reading, Massachusetts, 1981.

\item[\mbox{\rm[7]}] A. Schrijver,
Graph invariants in the spin model,
{\em Journal of Combinatorial Theory, Series B} 99 (2009) 502--511. 

\item[\mbox{\rm[8]}] B. Szegedy,
Edge coloring models and reflection positivity,
{\em Journal of the American Mathematical Society}
20 (2007) 969--988.

\end{itemize}
}

\end{document}